\documentclass[12pt]{amsart}

\usepackage{amsmath,amssymb,amsthm}
\newtheorem{Theorem}{Theorem}[section]

\newtheorem{Corollary}[Theorem]{Corollary}
\newtheorem{Proposition}[Theorem]{Proposition}
\newtheorem{Lemma}[Theorem]{Lemma}

\theoremstyle{definition}
\newtheorem{Example}[Theorem]{Example}
\newtheorem{Remark}[Theorem]{Remark}

\catcode`ð=\active \defð{\~{g}} \catcode`Ð=\active \defÐ{\~{G}}
\catcode`Ý=\active \defÝ{\. I} \catcode`ö=\active \defö{\"{o}}
\catcode`Ö=\active \defÖ{\"O} \catcode`ü=\active \defü{\"{u}}
\catcode`Ü=\active \defÜ{\"{U}} \catcode`Þ=\active \defÞ{\c{S}}
\catcode`þ=\active \defþ{\c{s}} \catcode`ý=\active \defý{{\i }}
\catcode`ç=\active \defç{\d{c}} \catcode`Ç=\active \defÇ{\d{C}}

\begin{document}

\title{Triviality of Symplectic $SU(2)$-actions on homology}

\author{Y{\i}ld{\i}ray Ozan}
\address{\hskip-\parindent
        Y{\i}ld{\i}ray Ozan\\
        Mathematics Department\\
        Middle East Technical University\\
        06531, Ankara\\
        Turkey}
\email{ozan@metu.edu.tr}


\thanks{The work presented in this paper has
been done during my stay at MSRI, Mathematical Sciences Research
Institute, Berkeley, in the Spring of 2004. The author is grateful
to the institute for the invitation and the kind hospitality. The
author would like to thank Dusa McDuff for pointing out some
errors in the earlier version of the paper and for several other
useful remarks.}

\begin{abstract}
Lalonde and McDuff showed that the natural action of the rational
homology of the group of Hamiltonian diffeomorphisms of a closed
symplectic manifold $(M, \omega)$ on the rational homology groups
$H_*(M,{\mathbb Q})$ is trivial. In this note, given a symplectic
action of $SU(2)$, $\phi:SU(2)\times M \rightarrow M$, we will
construct a symplectic fiber bundle $P_\phi\rightarrow {\mathbb
CP}^2$ with fiber $(M,\omega)$ and use it to construct the chains,
which bound the images of the homology cycles under the trace map
given by the $SU(2)$-action. It turns out that the natural chains
bounded by the $SU(2)$-orbits in $M$ are punctured ${\mathbb
CP}^2$'s, the counter parts of holomorphic discs bounding circles
in case of Hamiltonian circle actions. We will also define some
invariants of the action $\phi$ and do some explicit calculations.
\end{abstract}

\maketitle

\section{Introduction}
Let $\phi:G\times M \rightarrow M$ be a smooth action of a compact
Lie group $G$ on a smooth manifold $M$. The action induces a
homomorphism on homology, called the trace homomorphism,
$\partial_\phi:H_k(M,{\mathbb Q})\rightarrow H_{k+d}(M,{\mathbb
Q})$ defined as follows: If $\alpha \in H_{k}(M,{\mathbb Q})$ is a
class represented by a cycle $a:A\rightarrow M$, then
$\partial_\phi (\alpha)$ is the class in $H_{k+d}(M,{\mathbb Q})$
represented by the cycle $G\times A\rightarrow M$, $(g,x)\mapsto
\phi(g,a(x))$, where $d$ is the dimension of $G$. In general this
homomorphism is not trivial (just consider product spaces $G\times
M$). However, if $M$ is a closed symplectic manifold and $G$ is a
compact Lie group acting on $M$ in a Hamiltonian fashion, then it
is known that the homomorphism $\partial_\phi:H_k(M,{\mathbb
Q})\rightarrow H_{k+d}(M,{\mathbb Q})$ is trivial (cf. see
\cite{AB}). Later, Lalonde and McDuff have proved a stronger
result that the natural action of the homology of the group of
Hamiltonian diffeomorphisms of the closed manifold $(M, \omega)$
on the homology groups $H _*(M,{\mathbb Q})$,
$$H_k(\textmd{Ham}(M,\omega),{\mathbb Q})\times H_l(M,{\mathbb Q})
\rightarrow H_{k+l}(M,{\mathbb Q})$$ is trivial, for $k>0$
(\cite{LM,LMP}).

Below is the main result of this note, which determines the chains
bounded by the images of the trace homomorphism in the case of
$G=SU(2)$.

\begin{Theorem}\label{thm-A}
Let $\phi:SU(2)\times M \rightarrow M$ be a symplectic action on a
closed symplectic manifold $(M,\omega)$. Then there is a closed
symplectic manifold $(P_\phi,\omega_\phi)$,
which fibers over ${\mathbb CP}^2$ with fiber $M$ such that,
\begin{enumerate}
\item[i)] the rational homology of the fiber $M$ injects into the
rational homology of $P_\phi$, \item[ii)] the symplectic form
$\omega_\phi$ restricts to $\omega$ at each fiber, and
\item[iii)]
if $\alpha \in H_{k}(M,{\mathbb Q})$ is a class represented by a
cycle $a:A\rightarrow M$, then in the manifold $P_\phi$, the cycle
$$SU(2)\times A\rightarrow M\subseteq P_\phi$$ representing the class
$\partial_\phi (\alpha)$, bounds a chain of
the form ${\mathbb CP_0}^2 \times A\rightarrow P_\phi$, where
${\mathbb CP_0}^2={\mathbb CP}^2-Int(D^4)$ is the punctured
projective plane.
\end{enumerate}
In particular, the induced homomorphism on homology
$\partial_\phi:H_k(M,{\mathbb Q})\rightarrow H_{k+3}(M,{\mathbb
Q})$ is trivial.
\end{Theorem}

\begin{Remark}\label{rem-A}
Since $SU(2)$ is simply connected any symplectic $SU(2)$-action on
a symplectic manifold is Hamiltonian (\cite{C,MS}).
\end{Remark}

\begin{Example}\label{Example1}
Let $SU(2)$ act linearly on ${\mathbb CP}^2$ in the usual way (see
the next section). Blowing up the isolated fixed point of the
action we get an $SU(2)$-action on ${\mathbb CP}^2\sharp
\overline{\mathbb CP}^2$. The action is Hamiltonian since it is
algebraic. The orbit of a point with trivial stabilizer is a copy
of $SU(2)=S^3$, which separates the two copies of the projective
planes. So, the homology class represented by this orbit is
trivial and it bounds a punctured ${\mathbb CP}^2$, not a
$4$-ball.
\end{Example}

The next section is devoted to the proof of Theorem~\ref{thm-A}.
In the third section, we will construct some invariants of the
action $\phi:SU(2)\times M\rightarrow M$ and compute them in some
cases. Finally, we will mention some applications of these results
to the study of the topology of real algebraic varieties.

\section{Proof of Theorem~\ref{thm-A}}
Let $(M,\omega)$ be a closed $2n$-dimensional symplectic manifold
and $$\phi:SU(2)\times M\rightarrow M$$ a symplectic action. The
proof of Theorem~\ref{thm-A} consists of three parts. In the first
part, we will construct a smooth symplectic fiber bundle
$\pi_{S^4}:P^0_\phi\rightarrow S^4$ with fiber $(M,\omega)$ using
the action $\phi:SU(2)\times M\rightarrow M$ as the clutching
function. Moreover, the fibre bundle, both the total space and the
base, will have an $SU(2)$ and an $S^1$-action both preserving a
closed two form $\omega_{S^4}$ on $P^0_\phi$, which restricts to
$\omega$ on each fiber.  Moreover, the projection map will be
equivariant with respect to both actions.

In the second part, using a natural $SU(2)$ and $S^1$-equivariant
degree one map ${\mathbb CP}^2\rightarrow S^4$, where the $SU(2)$
and the $S^1$-action on the complex projective space are obtained
from the natural actions of these groups on ${\mathbb C}^2$, we
will pull back the bundle over the sphere to a bundle over
${\mathbb CP}^2$, which we will denote $\pi:P_\phi \rightarrow
{\mathbb CP}^2$. Using the pull back of the two form
$\omega_{S^4}$ on $P^0_\phi$ and the Fubuni-Study form on
${\mathbb CP}^2$ we will construct a symplectic form $\omega_\phi$
on $P_\phi$, which restricts to $\omega$ in each fiber. Moreover,
the $SU(2)$ and the $S^1-$actions on $P^0_\phi$ will induce
Hamiltonian actions on $P_\phi$.

In the last part, we will consider a symplectic reduction of the
total space of the bundle $\pi:P_\phi\rightarrow {\mathbb CP}^2$
using the $S^1$-action. Finally, Kirwan's Surjectivity Theorem on
symplectic quotients (\cite{K}) together with a topological
observation will finish the proof.

\subsection{Symplectic $M$-bundles over $S^4$ with structure
group $SU(2)$} \label{ss-sphere}

Given any smooth action $\phi:SU(2)\times M \rightarrow M$ we
define the smooth manifold $P^0_\phi$ as the identification space
$$D_+^4\times M \cup D_-^4\times M / (g,x)\sim (g,\phi(g,x)),
\hspace{0.2cm} \mbox{for} \ (g,x)\in\partial D_+^4\times M$$ where
we identify $\partial D_{\pm}^4$ with $S^3=SU(2)$. Note that we
have a fiber bundle $P^0_\phi \rightarrow S^4$ with fiber $M$,
induced by the projection maps $D_{\pm}^4 \times M \rightarrow
D_{\pm}^4$.

Another description for this bundle, which is more suitable to
define the $SU(2)$-action on, is as follows: Let ${\mathbb H}$
denote the quaternion line and $U({\mathbb H}^2)$ the set of unit
length vectors in ${\mathbb H}^2$. Also identify $SU(2)$ with the
set of unit quaternions ${\mathbb H}$. Let $L\rightarrow S^4$
denote the quaternion line bundle, whose unit disc bundle is the
$SU(2)$-bundle $U({\mathbb H}^2)\rightarrow S^4$,
$(v_1,v_2)\mapsto [v_1:v_2]$, for $(v_1,v_2)\in U({\mathbb H}^2)$.
Note that the latter map is nothing but the orbit map of the
$SU(2)$-action on $U({\mathbb H}^2)$ given by $(v_1,v_2)\mapsto
(v_1g,v_2g)$, for $g\in SU(2)$ and $(v_1,v_2)\in U({\mathbb
H}^2)$.

Similarly, on $U({\mathbb H}^2)\times M$ we have an $SU(2)$-action
defined by $$(g,(v_1,v_2),x)\mapsto ((v_1g,v_2g),\phi(g,x))$$ for
$g\in SU(2)$, $x \in M$ and $(v_1,v_2)\in U({\mathbb H}^2)$. The
projection map $U({\mathbb H}^2)\times M \rightarrow U({\mathbb
H}^2)$ is equivariant and taking quotients by the respective
actions on $U({\mathbb H}^2)\times M$ and $U({\mathbb H}^2)$ we
recover the fiber bundle $P^0_\phi \rightarrow S^4$.

There is a second $SU(2)$-action on $U({\mathbb H}^2)\times M$,
which commutes with the first one: $$(h,(v_1,v_2),x)\mapsto
((h^{-1}v_1,v_2),x)$$ for $h\in SU(2)$, $x \in M$ and
$(v_1,v_2)\in U({\mathbb H}^2)$. Clearly, this induces an action
on $S^4$ given as $$(h,[v_1:v_2])\mapsto [h^{-1}v_1:v_2]$$ which
makes the projection map equivariant.  Since the two actions
commute the second action induces actions on both the total space
and the base of the fiber bundle $P^0_\phi \rightarrow S^4$, which
makes the projection map equivariant.  Note that the action on
$S^4$ is free outside the poles, namely $[0:1]$ and $[1:0]$, the
only fixed points of the action.

There is also a left circle action on this space: Identify
${\mathbb H}$ with ${\mathbb C}^2$, on which $SU(2)$ acts by
matrix multiplication. Also regard $S^1$ as the set of matrices
$\{e^{i\theta}  I_2 \ | \ e^{i\theta}\in S^1\}$, where $I_2$ is
the $2\times 2$ identity matrix. Now let $S^1$ act on $U({\mathbb
H}^2)\times M$ by $$(e^{i\theta},(v_1,v_2),x)\mapsto
((e^{-i\theta} I_2v_1,v_2),x)$$ for $e^{i \theta}\in S^1$, $x \in
M$ and $(v_1,v_2)\in U({\mathbb H}^2)$. Since $e^{-i\theta} \ I_2$
commutes with all elements in $SU(2)$ we get an $S^1$-action on
both the total space and the base of the fiber bundle $P^0_\phi
\rightarrow S^4$, which makes the projection map equivariant. The
circle action on $S^4$ commutes with the $SU(2)$-action described
in the above paragraph and it is a free action outside the poles.

The Wang sequence for cohomology associated to the $M$-bundle
$P^0_\phi \rightarrow S^4$ yields the isomorphism
$$0=H^{-2}(M,{\mathbb Q})\rightarrow H^2(P^0_\phi,{\mathbb
Q})\stackrel{rest.}{\rightarrow} H^2(M,{\mathbb
Q})\stackrel{\partial_\phi^*}{\rightarrow} H^{-1}(M,{\mathbb
Q})=0$$ given by the restriction map, where the last map is the
dual of the trace homomorphism $\partial_\phi:H_k(M,{\mathbb
Q})\rightarrow H_{k+3}(M,{\mathbb Q})$ ($k=-1$ in this case). In
particular, there is unique cohomology class on $P^0_\phi$, which
restricts to the cohomology class $[\omega]$ on each fiber.
Indeed, we will construct a two form representing this cohomology
class, which we will use to get a symplectic form on the $M$-fiber
bundle over ${\mathbb CP}^2$.

\begin{Lemma}\label{lem-sphere}
There is an $SU(2)$ and $S^1$-invariant closed two form
$\omega_{S^4}$ on $P^0_{\phi}$, which restricts to $\omega$ on
each fiber, where the $SU(2)$ and $S^1$-actions are the ones
described in the above paragraphs. Moreover, the cohomology class
of $\omega_{S^4}$ is uniquely determined by these conditions.
\end{Lemma}

\begin{Remark}\label{rem-coupling}
The cohomology class $[\omega_{S^4}]$ is nothing but the usual
coupling class for the Hamiltonian fiber bundle
$P^0_{\phi}\rightarrow S^4$.
\end{Remark}

Before we prove this lemma we need some preliminaries. Let
$f:S^3\times M\rightarrow S^3\times M$ be the smooth map given by
$f(g,x)=(g,\phi(g,x))$, $(g,x)\in S^3\times M$. Consider the
following diagram, where $\pi_i$, $i=1,2$, are the projections
onto the second factors.
\begin{center}
$M$ \hspace{1.5cm} $M$

$\pi_1 \uparrow$ \hspace{1.7cm} $\uparrow \pi_2$

$S^3\times M  \ \stackrel{f}{\rightarrow} \ S^3\times M$
\end{center}
Using the decomposition $T_*(S^3\times M)=T_*S^3 \times T_*M$ we
will write any tangent vector $X$ on $S^3\times M$ as $X=(X_S,
X_M)$. Note that the differential of $f$ has the form $$f_*=\left(
\begin{array}{cc} Id & 0 \\ \frac{\partial \phi}{\partial g} &
\frac{\partial \phi}{\partial x} \end{array} \right).$$ Hence
$f_*((X_S,0))=(X_S, X_S^\sharp)$ and
$f_*((0,X_M))=(0,\phi_*(X_M))$, where $X_S^\sharp$ is the vector
field on $M$ generated by the vector $X_S$ via the action.

The $SU(2)=S^3$-action on $M$ is Hamiltonian means that there is a
smooth map $\mu :M\rightarrow su(2)^*$ such that for any vector
$X_S\in T_*S^3$ we have $$i_{X_S^\sharp}\omega=d(\mu(X_S)).$$
Since $SU(2)=S^3$ is parallelizable choosing a global frame $de_1,
\ de_2, \ de_3$ for the cotangent bundle for $S^3$ we can regard
$\mu$ as a one form on $S^3\times M$, namely
$$\mu(g,x)=A(x) \ de_1 + B(x) \ de_2+C(x) \ de_3$$ for $(g,x)\in
S^3\times M$.  One can easily check that
$d(\mu(X_S))=-i_{X_S}d\mu$.  Now we can state the next lemma.

\begin{Lemma}\label{lem-1}
$\pi_1^*(\omega)-(\pi_2\circ f)^*(\omega)$ is an exact two form on
$S^3\times M$, which vanishes on $T_*M$ identically.
\end{Lemma}

\begin{proof} Let $X=(X_S,X_M)$ and $Y=(Y_S,Y_M)$ be tangent vectors
at any point of $S^3\times M$. Then $\\ \mbox{\hspace{3.83cm}}
I=(\pi_1^*(\omega)-(\pi_2\circ f)^*(\omega)) \
((X_S,X_M),(Y_S,Y_M)) \\
\mbox{\hspace{4cm}} =\omega(X_M,Y_M)-\omega(X_S^\sharp
+\phi_*(X_M),Y_S^\sharp+\phi_*(Y_M))  \\
\mbox{\hspace{4cm}} = \omega(X_M,Y_M)-
\omega(\phi_*(X_M),\phi_*(Y_M))  \\
\mbox{\hspace{6cm}} -\omega(X_S^\sharp,Y_S^\sharp +\phi_*(Y_M))-
\omega(\phi_*(X_M),Y_S^\sharp)  \\
\mbox{\hspace{4cm}}  =-\omega(X_S^\sharp,Y_S^\sharp
+\phi_*(Y_M))-\omega(\phi_*(X_M),Y_S^\sharp), \\ $ because \
$\omega(\phi_*(X_M),\phi_*(Y_M))=\phi^*(\omega)(X_M,Y_M)=\omega(X_M,Y_M)$.
Note that this calculation already shows that
$\pi_1^*(\omega)-(\pi_2\circ f)^*(\omega)$ is identically zero on
$T_*M$.

Now using \ $i_{X_S^\sharp}\omega=d(\mu(X_S))$ \ we can write \\
$I=-d(\mu(X_S))(Y_S^\sharp +
\phi_*(Y_M))+d(\mu(Y_S))(\phi_*(X_M))$.
Since $$d(\mu(X_S))=-i_{X_S}d\mu$$  we get \\
$I=d\mu(X_S,Y_S^\sharp + \phi_*(Y_M))-d\mu(Y_S,\phi_*(X_M))
\\ \mbox{\hspace{4cm}} = d\mu(X_S,Y_S^\sharp)+d\mu(X_S,\phi_*(Y_M))+
d\mu(\phi_*(X_M),Y_S)$.

On the other hand, similar calculations yield \\
$d\mu(f_*(X),f_*(Y))=d\mu(X_S,Y_S^\sharp)+d\mu(X_S^\sharp,Y_S)
\\ \mbox{\hspace{4cm}} +d\mu(X_S,\phi_*(Y_M))+
d\mu(\phi_*(X_M),Y_S)$.

Comparing the two equations we deduce that \\
$I=d\mu(f_*(X),f_*(Y))-d\mu(X_S^\sharp,Y_S)$. For the last term we
can write \\
$d\mu(X_S^\sharp,Y_S)=-d\mu(Y_S,X_S^\sharp)=-(i_{Y_S}d\mu)(X_S^\sharp)
=d(\mu(Y_S))(X_S^\sharp) \\ \mbox{\hspace{4cm}} =(i_{Y_S^\sharp}
\omega)(X_S^\sharp)=\omega(Y_S^\sharp,X_S^\sharp)$.  So, we have
obtained  $$I=d\mu(f_*(X),f_*(Y))+\omega(X_S^\sharp,Y_S^\sharp).$$
Writing $\omega(X_S^\sharp,Y_S^\sharp)=I-(f^*(d\mu))(X,Y)$ we see
that the map $$(X,Y)\mapsto \omega(X_S^\sharp,Y_S^\sharp)$$ is a
closed two form on $S^3\times M$.  Moreover, it vanishes if $X_S$
or $Y_S$ is zero and hence it is identically zero on the $T_*M$
component of the tangent space. So, the de Rham cohomology class
represented by this closed two form evaluates zero on any two
dimensional homology class of the product $S^3\times M$ provided
that the homology class is represented by a cycle lying in some \
$\{pt\}\times M$. However, since $S^3$ has no first and second
homology, by the K{\"u}nneth formula, this de Rham class must be
trivial. Hence, there is a one form \ $u$ \ on $S^3\times M$ such
that $I=(f^*(d\mu))(X,Y)+du(X,Y)$. This finishes the proof of the
lemma.
\end{proof}

\begin{proof}[Proof of Lemma~\ref{lem-sphere}]
By the isomorphism obtained from the Wang sequence, the cohomology
class of $\omega_{S^4}$ is uniquely determined by that of $\omega$
(see the paragraph above Lemma~\ref{lem-sphere}).

We will regard the total space $P^0_\phi$ as the identification
space $${\mathbb R}_+^4\times M \ \cup \ {\mathbb R}_-^4\times M /
(t,g,x)\sim F(t,g,x)\doteq(t^{-1},g,\phi(g,x)),$$ for $(t,g,x)\in
({\mathbb R}_+^4-\{0\}) \times M$, where we identify ${\mathbb
R}^4-\{0\}$ with $(0,\infty)\times S^3$ in the obvious way. Let
$\pi_1$ and $\pi_2$ denote the projections onto the $M$ factors of
the products ${\mathbb R}_+^4\times M$ and ${\mathbb R}_-^4\times
M$, respectively.  Also, we will denote the projection of
$(0,\infty)\times S^3\times M$ onto the $S^3\times M$ component by
$\pi_{SM}$.

Let $\omega_i=\pi_i^*(\omega)$, \ $i=1,2$. Then by
Lemma~\ref{lem-1} $$\omega_1-F^*(\omega_2)=\pi_{SM}^*(dv_1)$$ for
some one form $v_1$ on $S^3\times M$. Let $v_2=(f^{-1})^*(v_1)$,
where $f$ is as in Lemma~\ref{lem-1}. So,
$\omega_1=F^*(\omega_2+dv_2)$ on $({\mathbb R}_+^4-\{0\}) \times
M$.  Let $\widetilde{\omega}_2=\omega_2+d(\rho(t)v_2)$, where
$\rho$ is a smooth function on ${\mathbb R}$, which vanishes on
$(-\infty,0.5]$ and equals one on $[0.75,\infty)$.

Now we have
$$F^*(\widetilde{\omega}_2)=F^*(\omega_2)+d(F^*(\rho(t)v_2))=
\omega_1 - dv_1 + d(\rho(1/t)v_1).$$  So, letting
$\widetilde{\omega}_1=\omega_1+d((\rho(1/t)-1)v_1)$ \ we obtain
$F^*(\widetilde{\omega}_2)=\widetilde{\omega}_1$.  By the choice
of the function $\rho$ the forms $\widetilde{\omega}_i$ are
defined on all of ${\mathbb R}^4\times M$, and indeed, they are
equal to $\omega_i$ on $D_{1/2}\times M$, where $D_{1/2}$ denotes
the disc of radius $1/2$ in ${\mathbb R}^4$. Moreover, since both
$dv_i$ and \ $dt$ \ vanishes on $T_*M$, each
$\widetilde{\omega}_i$ restricts to $\omega$ on each fiber
$\{pt\}\times M$. Hence, together they define a global closed two
form on $P^0_\phi$, say $\omega_{S^4}$, which restricts to
$\omega$ in each fiber.

The form $\omega_{S^4}$ may not be $SU(2)$-invariant. However, we
can average it over the $SU(2)$-orbits to get an $SU(2)$-invariant
form with the desired properties. Namely, let $dH$ denote the Haar
measure on $SU(2)$ with total volume one and define the average of
$\omega_{S^4}$ as the form $$(X,Y)\mapsto \int_{SU(2)}
(\psi^*(h,p) \ (\omega_{S^4}))\ (X,Y) \ dH$$ where
$\psi:SU(2)\times P^0_\phi \rightarrow P^0_\phi$ is the action map
and the integration is over $h\in SU(2)$ for fixed vectors $X,Y\in
T_p(P^0_\phi)$. Since on $U({\mathbb H}^2)\times M$ the action is
given by $(h,((v_1,v_2),x))\mapsto ((h^{-1}v_1,v_2),x)$ and the
restriction of $\omega_{S^4}$ to each fiber, which is a copy of
$M$, is just $\omega$, so will be the restriction of the average
of $\omega_{S^4}$. Once, the form is $SU(2)$-invariant then we can
average it to make also $S^1$-invariant in the same way. Since the
two actions commute, averaging over the $S^1$-orbits will not
spoil the $SU(2)-$invariance of the form. Also averaging commutes
with exterior derivative and therefore the averaged two form will
be still closed.
\end{proof}

\subsection{Symplectic $M$-bundles over ${\mathbb CP}^2$ with
structure group $SU(2)$}\label{ss-cp2}

By a fiber bundle $\pi: P_\phi \rightarrow {\mathbb CP}^2$, with
fiber $M$ and structure group $SU(2)$ we mean a group homomorphism
$\phi:SU(2)\rightarrow \textmd{Symp}(M,\omega)$, where the latter
is the group of symplectomorphisms of the symplectic manifold
$(M,\omega)$ and a principal $SU(2)$-bundle $P\rightarrow {\mathbb
CP}^2$ such that $P_\phi$ is obtained from $P$ via the
representation $\phi$ in the usual way.  The classifying space for
$SU(2)$-bundles is ${\mathbb HP}^{\infty}$, whose $7$th skeleton
is ${\mathbb HP}^1=S^4$ and therefore any principal $SU(2)$-bundle
over a closed $4$-manifold $N$ is obtained form the universal
bundle $L\rightarrow S^4$ by pulling it over $N$ by a map
$\xi:N\rightarrow S^4$.  Since $e(L)=c_2(L)\in H^4(S^4,{\mathbb
Z})$ is a generator we have $e(\xi^*(L))=c_2(\xi^*(L))=\deg(\xi)$.

Let $\xi:{\mathbb CP}^2\rightarrow S^4$ be the map given by the
formula
$$\xi([z_0:z_1:z_2])=(\frac{2\bar{z}_0z_1}{||z||^2},
\frac{2\bar{z}_0z_2}{||z||^2},\frac{|z_1|^2+|z_2|^2-|z_0|^2}{||z||^2})$$
where $||z||^2=|z_0|^2+|z_1|^2+|z_2|^2$, for $[z_0:z_1:z_2]\in
{\mathbb CP}^2$.  Here we consider the $4$-sphere as
$$S^4=\{(w_1,w_2,t)\in {\mathbb C}\times {\mathbb C} \times
{\mathbb R} \ | \ |w_1|^2+|w_2|^2+t^2=1 \}.$$  Note that
$\xi([0:z_1:z_2])=(0,0,1)$, the North pole and
$\xi([0:0:1])=(0,0,-1)$, the South pole. So, the map $\xi$ maps
the complex line $z_0=0$ to the North pole and is a diffeomorphism
onto its image outside the line $z_0=0$. In particular its degree
is one.

Consider the linear $SU(2)$-action on ${\mathbb CP}^2$ given as
$$[z_0:z_1:z_2]\mapsto [z_0:a'z_1+b'z_2:c'z_1+d'z_2]$$ where
$\left(\begin{array}{cc} a' & b' \\ c' & d' \end{array} \right)\in
SU(2)$ is the inverse of $\left(\begin{array}{cc} a & b
\\ c & d \end{array} \right)\in SU(2)$. Writing
$$w=(w_1,w_2)=(z_1/z_0,z_2/z_0)$$ we get
$$\xi([z_0:z_1:z_2])=(\frac{2w_1}{1+||w||^2},
\frac{2w_2}{1+||w||^2},1-\frac{2}{1+||w||^2}).$$ Hence, $\xi$
becomes equivariant if we endow $S^4$ with the $SU(2)$-action
given by $$(w_1,w_2,t)\mapsto (a'w_1+b'w_2,c'w_1+d'w_2,t).$$
However, the latter is just the action of Lemma~\ref{lem-sphere}.
Similarly, $\xi$ is $S^1$-equivariant where the $S^1$-action
${\mathbb CP}^2$ is given by
$$[z_0:z_1:z_2]\mapsto [z_0:e^{-i\theta}z_1:e^{-i\theta}z_2].$$

Let $P_\phi=\xi^*(P^0_\phi)$, the pull back of the $M$-bundle
$P^0_\phi\rightarrow S^4$ via the $SU(2)$ and the
$S^1$-equivariant map $\xi:{\mathbb CP}^2\rightarrow S^4$. Since
$\xi$ is an equivariant map the bundle $\pi:P_\phi\rightarrow
{\mathbb CP}^2$ gets both $SU(2)$ and $S^1$-actions, for which the
projection map $\pi$ is equivariant.  Moreover, the pull back
cohomology class $\xi^*(\omega_{S^4})$ is invariant with respect
to both actions and restricts to $\omega$ on each fiber.

Let $\omega_{FS}$ denote the Fubuni-Study symplectic form on
${\mathbb CP}^2$. The form $\pi^*({\omega_{FS}})$ is invariant
under the $SU(2)$ and the $S^1$-action on the complex projective
plane and is identically zero when restricted to each fiber
$\{pt\}\times M$. Hence for any positive large enough constant
$\kappa \gg 0$ the $2$-form
$\omega_\phi=\xi^*(\omega_{S^4})+\kappa \ \pi^*({\omega_{FS}})$ is
a symplectic form on $P_\phi$. Moreover, both actions on $P_\phi$
are Hamiltonian.  This is obvious for $SU(2)$ since it is simply
connected. For the $S^1$-action one can argue as follows: Since
averaging commutes with exterior derivative locally we have
$\omega_\phi=\pi_i^*(\omega)+dv +\kappa \ \pi^*({\omega_{FS}})$
for some equivariant one form \ $v$ \ on $P_\phi$.  Let $\chi$ be
a vector field generated by the $S^1$-action.  Since the form is
invariant the Lie derivative of $dv$ along $\chi$ will be zero.
Now by the Cartan formula we get
$i_{\chi^\sharp}dv=-d(i_{\chi^\sharp}v)$ and hence the
$S^1$-action on $P_\phi$ is also Hamiltonian.

\begin{Remark}\label{rem-L}

{\bf 1)} McDuff pointed out that the bundle $P_\phi \rightarrow
{\mathbb CP}^2$ can be also obtained as follows: Take a circle
subgroup $S^1$ of $SU(2)$ and consider the universal bundle
$M_{S^1}\rightarrow BS^1={\mathbb CP}^\infty$. Now the restriction
of this bundle to ${\mathbb CP}^2$ is the bundle $P_\phi
\rightarrow {\mathbb CP}^2$. To see this consider the fibration
$$S^2\rightarrow {\mathbb CP}^3=S^7/S^1\rightarrow S^7/SU(2)=
{\mathbb HP}^1=S^4.$$  The Gysin cohomology sequence for this
$S^2$-bundle yields that the restriction map ${\mathbb
CP}^2\rightarrow S^4$ has degree one.  Since
$P^0_\phi=S^7\times_{SU(2)}M$ and $P_\phi$ is the pullback of
$P^0_\phi$ via a degree one map ${\mathbb CP}^2\rightarrow S^4$
the assertion follows.

Indeed, more is true: Kedra and McDuff showed in \cite{KM} that a
homotopically trivial Hamiltonian circle action gives a nonzero
class in $\pi_4(B\textmd{Ham}(M,\omega))\otimes{\mathbb Q}$. The
class is defined as a Samelson product. Alternatively, one can see
this class as follows: The circle action gives a Hamiltonian
$(M,\omega)$-bundle over ${\mathbb CP}^2$ as described in above
paragraph. Since the action is homotopically trivial the bundle is
trivial over 2-skeleton ${\mathbb CP}^1$ of ${\mathbb CP}^2$. Now
gluing the trivial $(M,\omega)$-bundle over a three disc to the
bundle along a trivialization over ${\mathbb CP}^1$ we get an
$(M,\omega)$-bundle over $S^4$. Now the homotopy class of the
classifying map $S^4\rightarrow B\textmd{Ham}(M,\omega)$ of the
bundle over $S^4$ is the element found by Kedra and McDuff, whose
non triviality is proved in \cite{KM} using symplectic-Hamiltonian
characteristic classes.

{\bf 2)} By multiplying the last coordinate of the map
$\xi:{\mathbb CP}^2\rightarrow S^4$ by $-1$, if necessary, we can
arrange so that the pull back $SU(2)$-bundle $\xi^*(L)\rightarrow
{\mathbb CP}^2$ has $c_2=-1$. Since $SU(2)$-bundles are determined
by $c_2$ we see that $\xi^*(L)$ is smoothly isomorphic to
$O(1)\oplus O(-1)$. Therefore, the construction of $P_\phi$ could
be made just over ${\mathbb CP}^2$ using the bundle $O(1)\oplus
O(-1)$ with appropriate $SU(2)$ and $S^1$-actions.

{\bf 3)} We will orient the bundles $P^0_\phi$ and $P_\phi$ as
follows: The manifold $P_\phi$ is oriented with the orientation
coming from the symplectic form $\omega_\phi$. Since $P_\phi$ is
the pull back of $P^0_\phi$ via the map $\xi:{\mathbb CP}^2
\rightarrow S^4$, whose degree is chosen as above, the orientation
on $P_\phi$ induces one on $P^0_\phi$.
\end{Remark}

\subsection{Hamiltonians and symplectic reduction}
Let $\mu:P_\phi \rightarrow {\mathbb R}$ be a Hamiltonian for the
$S^1$-action on $P_\phi$. Recall that the $S^1$-equivariant map
$\xi:{\mathbb CP}^2\rightarrow S^4$ maps the line $z_0=0$ to the
North pole and sends the point $[0:0:1]$ to the South pole of the
sphere. So, over some small $S^1$-invariant disjoint tubular
neighborhoods $U$ and $V$ of the line $z_0=0$ and the point
$[0:0:1]$, respectively, the bundle $\pi:P_\phi\rightarrow {\mathbb
CP}^2$ is isomorphic to the product bundles $U\times M \rightarrow
U$ and $V\times M\rightarrow V$, where the $S^1$-action on the
$M$-factor is trivial. Moreover, by the construction, $\omega_\phi$
when restricted to $\pi^{-1}(U)$ and $\pi^{-1}(V)$, is just
$({\omega_\phi})_|=\omega+\kappa \ \pi^*(\omega_{FS})$. Since the
action on the $M$-factor is trivial we see that the moment map
restricted to $\pi^{-1}(U)$ and $\pi^{-1}(V)$ is just a multiple of
the moment map $\mu_0:{\mathbb CP}^2\rightarrow {\mathbb R}$ of the
$S^1$-action on ${\mathbb CP}^2$ plus a constant, which depends only
on the open set $U$ or $V$; i.e., $\mu(p)=\kappa \
\pi(\mu_0(p))+C(p)$, for all $p\in \pi^{-1}(U) \cup \pi^{-1}(V)$,
where $C$ is a locally constant function on the union
$\pi^{-1}(U)\cup \pi^{-1}(V)$.

We are ready now to prove the main theorem.
\begin{proof}[Proof of Theorem~\ref{thm-A}]
Replacing the Hamiltonians by adding constants if necessary we can
assume that \ $0$ \ is a regular value for $\mu$ and hence for
$\mu_0$ such that $\mu_0^{-1}(0)=S^3$ lies in $U$. Note that this
$S^3$ divides ${\mathbb CP}^2$ into two pieces. By multiplying all
the symplectic forms with \ $-1$ \ if necessary we can assume that
${\mathbb CP}^2_0=\mu_0^{-1}((-\infty,0])$ is a closed tubular
neighborhood of the line $z_0=0$ and $D^4=\mu_0^{-1}([0,\infty))$
is a closed $4$-ball with common boundary $S^3=\mu^{-1}(0)$. The
$M$-fiber bundle over these pieces are just products and
$$P_\phi={\mathbb CP}^2_0\times M \cup D^4\times M / (g,x)\sim (g,\phi(g,x)),
\hspace{0.2cm} \mbox{for} \ (g,x)\in\partial ({\mathbb
CP}^2_0)\times M.$$ Let $\alpha \in H_{k}(M,{\mathbb Q})$ be a
class represented by a cycle $a:A\rightarrow M$.  We need to show
that the class $\partial_\phi (\alpha)$, represented by the cycle
$S^3\times A\rightarrow M$, $(g,x)\mapsto \phi(g,a(x))$, is
trivial in $H_{k+3}(M,{\mathbb Q})$. We can clearly view
$S^3\times A$ as a subset of the boundary of ${\mathbb
CP}^2_0\times M$. The identification in the above decomposition of
$P_\phi$ maps $S^3\times A$ into the other piece by the map
$(g,x)\mapsto \phi(g,a(x))$. On the other hand, the radial
contraction of $D^4$ to its center \ $\{ 0 \}$ \ induces a radial
contraction of $D^4\times M$ to $\{ 0 \}\times M$. Moreover, the
composition of the identification map with the contraction will
map $S^3\times A$ into $M$ exactly via the map $(g,x)\mapsto
\phi(g,a(x))$. Since $S^3\times A=\partial ({\mathbb CP}^2_0\times
A)$ the class $\partial_\phi (\alpha)$ is trivial in
$H_{k+3}(P_\phi,{\mathbb Q})$.

Now consider the symplectic quotient $\mu^{-1}(0)/S^1$, which is
equal to the product $S^3/S^1\times M=S^2\times M$, because $S^1$
acts trivially on $M$ by the construction of the $S^1$-action.  By
the Kirwan's Surjectivity Theorem (\cite{K}) the map, induced by the
inclusion $\mu^{-1}(0)\subseteq P_\phi$,
$\mathcal{K}:H_{S^1}^i(P_\phi,{\mathbb Q})\rightarrow H^i(S^2\times
M,{\mathbb Q})$ is onto, for all $i$. So the restriction map
$H_{S^1}^i(P_\phi,{\mathbb Q})\rightarrow H^i(M,{\mathbb Q})$ is
surjective. Hence the map in homology induced by the inclusion of a
fiber $H_i(M,{\mathbb Q})\rightarrow H^{S^1}_i(P_\phi,{\mathbb Q})$
is injective. By the above paragraph $\partial_\phi (\alpha)$ is
trivial in $H_{k+3}(P_\phi,{\mathbb Q})$. However, both the cycle
$S^3\times A$ and the chain bounding it, ${\mathbb CP}^2_0\times A$,
are $S^1$-equivariant.  Hence the class $\partial_\phi (\alpha)$ is
also trivial in $H^{S^1}_{k+3}(P_\phi,{\mathbb Q})$. This implies
that $\partial_\phi (\alpha)$ must vanish in $H_{k+3}(M,{\mathbb
Q})$.
\end{proof}

\section{Some invariants of the SU(2)-action}
In this section we will study the sections of the bundles
$P^0_\phi$ and $P_\phi$, define some invariants of the
$SU(2)$-action on $(M,\omega)$, make some computations and mention
some applications to the study of the topology of real algebraic
varieties.

The orientations on the manifolds $P^0_\phi$ and $P_\phi$, which
we will need when we consider integrals over them, are the ones
described in Remark~\ref{rem-L}.

\subsection{Sections of $P^0_\phi$}
The lemma below describes the $SU(2)$-equivariant (with respect to
the $SU(2)$-action on $P^0_\phi$ and on $S^4$ described in
Subsection~\ref{ss-sphere}) sections of the bundle
$P^0_\phi\rightarrow S^4$ up to homotopy.

\begin{Lemma}\label{lem-Ssection}
There is an $SU(2)$-equivariant section $s:S^4\rightarrow P^0_\phi$
if and only if the $SU(2)$-action on $M$ has a fixed point.
Moreover, if $s_i:S^4\rightarrow P^0_\phi$, $i=1,2$, are any two
sections (not necessarily equivariant) then the difference of
$(s_2)_*([S^4]) -(s_1)_*([S^4])$ as a homology class is in the image
of the composition map $\pi_4(M)\rightarrow
\pi_4(P^0_\phi)\rightarrow H_4(P^0_\phi,{\mathbb Z})$, induced by
the inclusion of a fiber.
\end{Lemma}

\begin{proof}
Let $l(t):[-1,1]\rightarrow S^4$ be a one to one geodesic arc from
the point $(0,0,-1)$ to the point $(0,0,1)$. If $s:S^4\rightarrow
P^0_\phi$ is an equivariant section then $s$ is determined
completely by its values $s(l(t))$, $t\in [-1,1]$. On the other
hand, the points $(0,0,\pm 1)$ are the fixed points of the action
on the sphere and hence the points $s(0,0,\pm 1)$ are in the fixed
point set of action on $M$.  Moreover, since the action on $S^4$
is free outside the poles, any section $s$ defined on the arc
$l(t)$ with $s(0,0,\pm 1)\in M$ fixed points of the
$SU(2)$-action, extends uniquely to a section. Indeed, the section
$s:S^4\rightarrow P^0_\phi$ is just the trace of the section
$s(l(t))$, $t\in [-1,1]$, under the $SU(2)$-action on $P^0_\phi$.

The second statement follows the long exact sequence corresponding
to the fibration $M\rightarrow P^0_\phi \rightarrow S^4$,
$$\cdots \rightarrow \pi_4(M)\rightarrow \pi_4(P^0_\phi) \rightarrow
\pi_4(S^4)={\mathbb Z} \rightarrow \cdots .$$
\end{proof}

Theorem~\ref{thm-A} implies the following result.

\begin{Proposition}\label{prop-Ssection}
If the $SU(2)$ action on $M$ is symplectic then the homology class
$s_*([S^4])$ of an equivariant section $s:S^4\rightarrow P^0_\phi$
is determined only by the connected components of the fixed point
set containing the fixed points $s(0,0,-1)$ and $s(0,0,1)$.
\end{Proposition}

\begin{proof}
Assume the set up in the proof of the Lemma~\ref{lem-Ssection}. If
$s_1$ and $s_2$ are two such sections with
$s_1(0,0,-1)=s_2(0,0,-1)$ and $s_1(0,0,1)=s_2(0,0,1)$ then the
difference homology class can be identified with the trace of a
loop in $M$ based at one of these two fixed points.  However, by
Theorem~\ref{thm-A} the latter is trivial. Now assume that
$s_1(0,0,-1)$, $s_2(0,0,-1)\in F_0$ and $s_1(0,0,1)$,
$s_2(0,0,1)\in F_1$, for some connected components $F_0$ and $F_1$
of the fixed point set. Join the fixed points in $F_0$ and $F_1$
by arcs contained completely in the fixed point sets. The trace of
a path that lies in the fixed point set is just the path itself
and hence these arcs do not contribute to the difference of the
homology classes.  This finishes the proof.
\end{proof}

We will call a cohomology class $u\in H^4(M,{\mathbb Q})$ {\it
monotone} if it vanishes on spherical classes, i.e., on the image
of $\pi_4(M)\rightarrow H_4(M,{\mathbb Q})$.

Following \cite{LMP,S}, we denote the Chern classes of the
vertical tangent bundle
$$T_*^{vert}=\ker(\pi_*:T_*P^0_\phi\rightarrow T_*S^4)$$ by
$c^\phi_i$. These classes are clearly invariants of the action
$\phi$ and hence, so is any integral of the form
$$I^0(k,k_1,\cdots,k_n)=\int_{P^0_\phi} \omega_{S^4}^k \ (c_1^\phi)^{k_1}
\cdots (c_n^\phi)^{k_n}$$ where $k$, $k_i$ are non negative
integers with $2k+2k_1+\cdots +2nk_n=2n+4=\dim(P^0_\phi)$.

Let $x\in M$ be any fixed point of the $SU(2)$-action. Then the
function $v\mapsto (v,x)$, $v\in S^4$, defines a section, say
$s_x:S^4\rightarrow P^0_\phi$. Note that the pull back bundle over
$S^4$ of the vertical bundle via the section $s_x$ is nothing but
the associated complex vector bundle of the principal
$SU(2)$-bundle $U({\mathbb H}^2)\rightarrow S^4$ (see
Subsection~\ref{ss-sphere}) corresponding to the representation of
$SU(2)$ on tangent space $T_xM$. We have then the following result
about the representations of $SU(2)$ on tangent spaces of the
fixed points, which follows easily from Lemma~\ref{lem-Ssection}.

\begin{Corollary}\label{cor-sphere}
Let $(M,\omega)$ and $\phi$ be as above. Assume that $c_2(M)$ is a
monotone class. Then for any two fixed points $x_1$ and $x_2$ of
the action on $M$, the vector bundles over $S^4$, corresponding to
the $SU(2)$-representations at the tangent spaces $T_{x_i}$,
$i=1,2$, have the same second Chern class.
\end{Corollary}

\begin{Remark}\label{rem-sp}
Lemma~\ref{lem-Ssection} and the above corollary are valid indeed
for any smooth action on $M$ since we do not make use of the
symplectic form at all.
\end{Remark}

\begin{Example}\label{Example2}
{\bf 1)} If the action $\phi:SU(2)\times M\rightarrow M$ is
trivial then the integrals $$I^0(k,k_1,\cdots,k_n)=\int_{P^0_\phi}
\omega_{S^4}^k \ (c_1^\phi)^{k_1} \cdots (c_n^\phi)^{k_n}$$ are
all zero, because in this case $P^0_\phi=S^4\times M$ and all the
forms in the integral are trivial on $T_*S^4$.

{\bf 2)} Consider the standard action of $SU(2)$ on $({\mathbb
CP}^2,\omega)$, where $\omega=\omega_{FS}$, the Fubuni-Study
metric. Note that $c_2(\mathbb CP^2)$ is a monotone class because
$\pi_4({\mathbb CP}^2)=0$. It follows from the Wang sequence of
the fibration that $c_2^\phi=\lambda
[\omega_{S^4}]^2+\pi_{S^4}^*(v)$ for some $v\in H^4(S^4,{\mathbb
Q})$ and real number $\lambda$. Let $x_0=[1:0:0]$, the only fixed
point of the action, and denote the corresponding section of
$P^0_\phi\rightarrow S^4$ by $s_{x_0}$. Then
$$c_2^\phi([s_{x_0}])=\lambda [\omega_{S^4}]^2([s_{x_0}])+
\pi_{S^4}^*(v)([s_{x_0}])=\pi_{S^4}^*(v)([s_{x_0}])=v([S^4]),$$
where the second equality follows from the fact
$s_{x_0}^*([\omega_{S^4}])=0$ and the third one from that
$\pi_{S^4}\circ s_{x_0}=id_{S^4}$. It follows from the
representation of $SU(2)$ on the tangent space $T_{x_0}{\mathbb
CP}^2$ that $c_2^\phi([s_{x_0}])=c_2(L)=-1$, where $L\rightarrow
S^4$ is the canonical $SU(2)$-bundle (see
Section~\ref{ss-sphere}).

Now, let us calculate the invariants $I^0(k,k_1,k_2)$ for this
action.  Theorem~\ref{thm-A} and the Wang sequence implies that
$H_4(P_\phi^0, {\mathbb Q})\simeq H_4(S^4, {\mathbb Q})\oplus
H_4({\mathbb CP}^2, {\mathbb Q})$, whose generators are the
section class $[s_{x_0}]$ and the fiber class $[{\mathbb CP}^2]$,
respectively. Note that the normal bundle to the section $s_{x_0}$
in $P_\phi^0$ is just the vertical bundle and hence the Euler
class of the normal bundle is the restriction of $c_2^\phi$ to
$s_{x_0}$. In particular, the self intersection of $[s_{x_0}]$ is
$-1$.  On the other hand, clearly $[{\mathbb CP}^2]\cdot [{\mathbb
CP}^2]=0$ and $[{\mathbb CP}^2]\cdot [s_{x_0}]=1$.

Let $\alpha\in H_4(P_\phi^0,{\mathbb Q})$ denote the Poincar\`{e}
dual of $c_2^\phi$. Since $c_2^\phi([{\mathbb CP}^2])=3$ and
$c_2^\phi([s_{x_0}])=-1$ we see that $\alpha= 3[s_{x_0}]+
2[{\mathbb CP}^2]$. Hence
$$I^0(0,0,2)=\int_{P_\phi^0}(c_2^{\phi})^2= \alpha \cdot \alpha =3.$$

Recall the Wang sequence for cohomology associated to the
$M$-bundle $P^0_\phi \rightarrow S^4$, which yields the
isomorphism
$$0=H^{-2}(M,{\mathbb Q})\rightarrow H^2(P^0_\phi,{\mathbb
Q})\stackrel{rest.}{\rightarrow} H^2(M,{\mathbb
Q})\stackrel{\partial_\phi^*}{\rightarrow} H^{-1}(M,{\mathbb
Q})=0.$$ Hence any cohomology class in $H^2(M,{\mathbb Q})$
extends to a class in $H^2(P^0_\phi,{\mathbb Q})$ uniquely. In
particular, since $[\omega_{FS}]=c_1({\mathbb CP}^2)$ we see that
$[\omega_{S^4}]=c_1^\phi$. Clearly,
$(c_1^\phi)^2([s_{x_0}])=(s_{x_0}^*(c_1^\phi))^2=0$ and
$(c_1^\phi)^2([{\mathbb CP}^2])=9$. So, if $\beta$ denotes the
Poincar\`{e} dual of $(c_1^\phi)^2=[\omega_{S^4}]^2$ then
$\beta=9[s_{x_0}]+9[{\mathbb CP}^2]$.  Therefore, for nonnegative
integers $i+j=2$
$$I^0(i,j,1)=I^0(2,0,1)=\int_{P_\phi^0}\omega_{S^4}^2 (c_2^{\phi})=
\alpha \cdot \beta =18$$ and for nonnegative integers $i+j=4$
$$I^0(i,j,0)=I^0(4,0,0)=\int_{P_\phi^0}\omega_{S^4}^4= \beta \cdot
\beta =81.$$
\end{Example}

Note that we can define analogous integrals over $P_\phi$: There
is a unique cohomology class $u_\phi=\xi^*([\omega_{S^4}])\in
H^2(P_\phi,{\mathbb Q})$, which restricts to $[\omega]$ on each
fiber such that the integration of $u_\phi^{n+1}$ along the fiber
is zero. Similarly, we define
$$I(k,k_1,\cdots,k_n)=\int_{P_\phi} u_\phi^k \ (c_1^\phi)^{k_1}
\cdots (c_n^\phi)^{k_n}$$ where $k$, $k_i$ are non negative
integers with $2k+2k_1+\cdots +2nk_n=2n+4=\dim(P_\phi)$.  Note
that the two invariants are indeed equal, where the first one may
be more suitable for computations. However, if the action has no
fixed points then the bundle over $S^4$ may have no section, as
the next example shows.

\begin{Example}\label{exam-nosection}
The bundle corresponding to the standard action of $SU(2)$ on
${\mathbb CP}^1$, $P^0_\phi\rightarrow S^4$, has no section.
Otherwise, by deforming the section to a constant, say $[1:0]$,
over $D^4$ we would arrive at the contradiction that the map
$SU(2)\rightarrow {\mathbb CP}^1$ given by
$$g=\left(\begin{array}{cc} a & b
\\ c & d \end{array} \right)\mapsto
[\bar{a}:\bar{b}]=\phi(g,[1:0])$$ is homotopically trivial, which
is nothing but basically the Hopf map (see Section~\ref{ss-cp2}).

On the other hand, as we will see in the next section that the
bundle $\pi:P_\phi\rightarrow {\mathbb CP}^2$ has always a
section.
\end{Example}

\subsection{Sections of $P_\phi$}
Since the $SU(2)$ action on ${\mathbb CP}^2$ has a fixed point, an
equivariant section exists if and only if the $SU(2)$-action on
$M$ has a fixed point. Note also that, any equivariant section of
$P^0_\phi\rightarrow S^4$ pulls back to an equivariant section of
$P_\phi \rightarrow {\mathbb CP}^2$. Indeed, these pull back
sections are those equivariant sections $s:{\mathbb
CP}^2\rightarrow P_\phi$ such that $s(z_0=0)$ and $s([1:0:0])$ are
two fixed points of the action on $M$ (since the map $\xi:{\mathbb
CP}^2\rightarrow S^4$ maps the line $z_0=0$ to a pole of $S^4$ the
bundle $P_\phi$ restricted to the line $z_0=0$ is trivial). In
particular, for this class of equivariant sections of $P_\phi$ the
analogues of Lemma~\ref{lem-Ssection},
Proposition~\ref{prop-Ssection} and Corollary~\ref{cor-sphere}
will hold.

Another source of equivariant sections of the bundle is the set of
points of $M$ whose stabilizers is a circle. Namely, let
$H=\textmd{Stab}_{SU(2)}([0:1:1])$ and consider the path
$r(t)=[1-t:t:t]$, $t\in [0,1]$, in ${\mathbb CP}^2$. Let $x_0$,
$x_1$ be points in $M$ with $\textmd{Stab}_{SU(2)}(x_0)=SU(2)$
($x_0$ is a fixed point) and
$H\leqslant\textmd{Stab}_{SU(2)}(x_1)$. Choose a section of the
bundle over the path $r(t)$ with $s([1:0:0])=x_0$ and
$s([0:1:1])=x_1$. Then this extends uniquely to an equivariant
section $s:{\mathbb CP}^2\rightarrow P_\phi$. Moreover, any
equivariant section is of this form and the analogues of
Lemma~\ref{lem-Ssection} and Corollary~\ref{cor-sphere} will hold
in this case also.

Unlike the bundle $P^0_\phi\rightarrow S^4$ the bundle over
${\mathbb CP}^2$ has always a section.

\begin{Lemma}\label{lem-section}
Let $p_0$ be any point in the fiber $\pi^{-1}([1:0:0])=M$. Then
the bundle $\pi:P_\phi\rightarrow {\mathbb CP}^2$ has a section
$s:{\mathbb CP}^2\rightarrow P_\phi$ with $s([1:0:0])=p_0$.
\end{Lemma}

\begin{proof}
Recall the decomposition of $P_\phi$ from the proof of the
Theorem~\ref{thm-A}
$$P_\phi={\mathbb CP}^2_0\times M \cup D^4\times M /
(g,x)\sim (g,\phi(g,x)),\hspace{0.2cm} \mbox{for} \
(g,x)\in\partial ({\mathbb CP}^2_0)\times M.$$ We define a section
$s:{\mathbb CP}^2\rightarrow P_\phi$ as follows: For $v\in D^4$
let $s(v)=(v,p_0)$.  Over $\partial(D^4)$ the section looks like
$v\mapsto (v,p_0)$ and over $\partial({\mathbb CP}^2_0)$ it is
given by the formula $v\mapsto (v,\phi(v^{-1},p_0))$.

Since the action of the maximal torus
$H=\textmd{Stab}_{SU(2)}([0:1:1])$ on $M$ is also Hamiltonian it
has a fixed point, say $p_1\in M$. Hence,
$\textmd{Stab}_{SU(2)}(p_1)$ contains $H$.  Let
$\sigma:[0,1]\rightarrow M$ be a path from $p_0$ to $p_1$.  The
$SU(2)$-orbit of this path,
$$SU(2)\times [0,1]\rightarrow M, \ (g,t)\mapsto \phi(g,\sigma(t))$$ gives a map
$\beta:{\mathbb CP}^2_0\rightarrow M$, whose restriction to the
boundary $\partial ({\mathbb CP}^2_0)$ is the orbit of $p_0$. Now
we can extend the section over ${\mathbb CP}^2_0$ as $v\mapsto
(v,\beta(v))$.
\end{proof}

\begin{Remark}\label{rem-existence}
Note that if the point $p_0\in M$ is not a fixed point of the
$SU(2)$-action then the section of the above lemma is not
equivariant.  In particular, if $M={\mathbb CP}^1$ with the
$SU(2)$-action as in Example~\ref{exam-nosection}, then neither
bundles have an equivariant section.
\end{Remark}

We believe that $J$-holomorphic sections of $P_\phi\rightarrow
{\mathbb CP}^2$ deserve some attention also.

\subsection{Algebraic actions on real algebraic varieties}\label{ss-real}
The result mentioned in the introduction that the natural action
of the homology of the group of Hamiltonian diffeomorphisms of a
closed symplectic manifold $(M, \omega)$ on the homology groups
$H_*(M,{\mathbb Q})$,
$$H_k(\textmd{Ham}(M,\omega),{\mathbb Q})\times H_l(M,{\mathbb Q})
\rightarrow H_{k+l}(M,{\mathbb Q})$$ is trivial (\cite{LM,LMP})
has an immediate consequence in the study of topology of real
algebraic varieties: Let $X$ be a nonsingular compact real
algebraic variety with a nonsingular projective complexification
$i:X \rightarrow X_{\mathbb C}$. Clearly $X_{\mathbb C}$ carries a
K\"ahler and hence a symplectic structure such that $X$ is a
Lagrangian submanifold. Define $KH_i(X,{\mathbb Q})$ as the kernel
of the homomorphism
$i_*:H_i(X,{\mathbb Q}) \rightarrow H_i(X_{\mathbb C},{\mathbb
Q})$ and $ImH^i(X,{\mathbb R})$ as the image of the
homomorphism
$i^*:H^i(X_{\mathbb C},{\mathbb Q}) \rightarrow H^i(X,{\mathbb
Q})$. In \cite{O1,O2} it is shown that both $KH_i(X,{\mathbb Q})$
and $ImH^i(X,{\mathbb Q})$ are independent of the projective
complexification $i:X \rightarrow X_{\mathbb C}$ and thus (entire
rational) isomorphism invariants of $X$. We know also that the
natural linear action of a unitary group on a complex projective
variety is Hamiltonian. We have then the following corollary.

\begin{Corollary}\label{cor-real}
Let $X$ and $X_{\mathbb C}$ be as above and $G$ be a compact Lie
group acting unitarily on $X_{\mathbb C}$, leaving the real part
$X$ invariant.  Then the image of the trace map
$$H_k(G,{\mathbb Q})\times H_l(X,{\mathbb Q}) \rightarrow
H_{k+l}(X,{\mathbb Q})$$ lies in $KH_{k+l}(X,{\mathbb Q})$.
\end{Corollary}

\providecommand{\bysame}{\leavevmode\hboxto3em{\hrulefill}\thinspace}

\end{document}